\documentclass[conference]{IEEEtran}
\IEEEoverridecommandlockouts

\usepackage{cite}
\usepackage[]{hyperref}
\usepackage{mathptmx, enumerate}
\usepackage{amsmath,amssymb,amsthm,amsbsy} %
\usepackage[dvipdfmx]{graphicx}
\usepackage{epstopdf}
\usepackage{mathrsfs}
\usepackage{mathtools}
\usepackage{bm}
\usepackage{yhmath}
\usepackage{color}
\usepackage{arydshln}
\usepackage[hang,small]{caption}
\usepackage[subrefformat=parens]{subcaption}
\usepackage{comment}
\mathtoolsset{showonlyrefs=true}

\newtheorem{theorem}{Theorem}

\newtheorem{proposition}{Proposition}

\theoremstyle{definition}

\newtheorem{example}{Example}
\newtheorem{fact}{Fact}
\newtheorem{problem}{Problem}
\newtheorem{remark}{Remark}
\newtheorem{assumption}{Assumption}

\newcommand{\setN}{\mathbb{N}}
\newcommand{\setR}{\mathbb{R}}
\newcommand{\setPR}{\mathbb{R}_{++}}
\newcommand{\setNNR}{\mathbb{R}_{+}}

\newcommand{\ip}[2]{\left\langle #1 , #2 \right\rangle}
\newcommand{\norm}[1]{\left\| #1 \right\|}

\newcommand{\ran}{\operatorname{ran}}

\newcommand{\nullsp}{\operatorname{null}}
\newcommand{\setLO}[2]{\mathcal{B}\left( #1 , #2 \right)}

\newcommand{\Id}{\operatorname{Id}}
\newcommand{\zeroMatrix}{\mathrm{O}}

\newcommand{\ri}{\operatorname{ri}}

\DeclareFontFamily{U}{mathx}{\hyphenchar\font45}
\DeclareFontShape{U}{mathx}{m}{n}{
      <5> <6> <7> <8> <9> <10>
      <10.95> <12> <14.4> <17.28> <20.74> <24.88>
      mathx10
      }{}
\DeclareSymbolFont{mathx}{U}{mathx}{m}{n}
\DeclareMathSymbol{\bigtimes}{1}{mathx}{"91}

\newcommand{\dom}{\operatorname{dom}}

\newcommand{\prox}{\operatorname{Prox}}
\newcommand{\fix}{\operatorname{Fix}}

\newcommand{\spH}{\mathcal{H}}
\newcommand{\spK}{\mathcal{K}}

\newcommand{\spX}{\mathcal{X}}
\newcommand{\spY}{\mathcal{Y}}
\newcommand{\spZ}{\mathcal{Z}}
\newcommand{\sptildeZ}{\widetilde{\mathcal{Z}}}

\newcommand{\spfrakZ}{\mathfrak{Z}}

\newcommand{\opP}{\mathfrak{P}}
\newcommand{\opL}{\mathfrak{L}}

\newcommand{\opC}{\mathfrak{C}}

\newcommand{\setS}{\mathcal{S}}
\newcommand{\minimize}[0]{\operatornamewithlimits{minimize\ }}
\newcommand{\argmin}{\operatornamewithlimits{argmin\ }}

\newcommand{\xstar}{{x^\star}}

\DeclareMathAlphabet\mathbfcal{OMS}{cmsy}{b}{n}

\renewcommand{\baselinestretch}{0.94}

\begin{document}

\title{A Convexity Preserving Nonconvex Regularization\\ for Inverse Problems under Non-Gaussian Noise
\thanks{This work was supported by JSPS Grants-in-Aid (19H04134, 23KJ0945, 24K23885).}
}

\author{\IEEEauthorblockN{Wataru Yata, Keita Kume, Isao Yamada}
\IEEEauthorblockA{\textit{Dept. of Information and Communications Engineering} \\
\textit{Institute of Science Tokyo}, Tokyo, Japan \\
Email: \{yata, kume, isao\}@sp.ict.e.titech.ac.jp}
}

\maketitle
\begin{abstract}
    We propose a nonconvexly regularized convex model for linear regression problems under non-Gaussian noise. The cost function of the proposed model is designed with a possibly non-quadratic data fidelity term and a nonconvex regularizer via the generalized Moreau enhancement of a seed convex regularizer. We present sufficient conditions (i) for the cost function of the proposed model to be convex over the entire space, and (ii) for the existence of a minimizer of the proposed model. 
    Under such conditions, we propose a proximal splitting type algorithm with guaranteed convergence to a global minimizer of the proposed model. 
    As an application, we enhance nonconvexly a convex sparsity-promoting regularizer in a scenario of simultaneous declipping and denoising.
\end{abstract}

\begin{IEEEkeywords}
generalized Moreau enhancement, convex optimization, non-quadratic data fidelity, proximal splitting.
\end{IEEEkeywords}
\section{Introduction} 
Sparsity-aware estimation of a target signal $x^\star\in\spX$ from its noisy observation:
\begin{align}
  \label{eq:linear-regression}
  y = Ax^\star +\varepsilon\in\spY
\end{align} 
is a central goal in inverse problems and signal processing \cite{elad2010,theodoridis2024}, where $\spX$ and $\spY$ are finite-dimensional real Hilbert spaces, $A\in\setLO{\spX}{\spY}$ is a known linear operator and $\varepsilon\in\spY$ is noise (see Section \ref{sec:notation} for mathematical notation).

In this paper, we consider the following nonconvexly regularized convex (NRC) model\footnote{The cost function $J$ in \eqref{eq:cLiGME-w-general-fidelity} enjoys the convexity by a technical condition (see Proposition \ref{prop:overall-convexity}).} for estimation of $\xstar$ in \eqref{eq:linear-regression} by using its prior knowledge on the sparsity of $\opL \xstar$ with a certain known linear operator $\opL$.
\begin{problem}[Target NRC model]
  \label{prob:cLiGME-w-general-fidelity}
  Let $\spX,\spY,\spZ,\sptildeZ$ and $\spfrakZ$ be finite-dimensional real Hilbert spaces.
  Under the following settings (see Section \ref{sec:notation} for definitions of \textit{italicized technical terms}):
  \begin{enumerate}[(a)]
    \item  $\mathbf{C} (\subset \spfrakZ)$ is a \textit{simple} closed convex set and $\opC\in\setLO{\spX}{\spfrakZ}$;
    \item $A\in\setLO{\spX}{\spY}$ and $f\in\Gamma_0(\spY)$ is continuously differentiable over $\spY$;
    \item We use $\mu\in\setPR$, $\opL\in\setLO{\spX}{\spZ}$ and $B\in\setLO{\spZ}{\sptildeZ}$. 
    $\Psi\in\Gamma_0(\spZ)$ is \textit{coercive} and \textit{prox-friendly};
    \item $\dom(\Psi\circ\opL)\cap \opC^{-1}(\mathbf{C})\neq\emptyset$,
  \end{enumerate}
  consider a nonconvexly regularized convex model:
  \begin{equation}
    \label{eq:cLiGME-w-general-fidelity}
    \minimize_{\opC x \in\mathbf{C}} J(x) := f\circ A(x) + \mu \Psi_B\circ\opL (x),
  \end{equation}
  where $\Psi_B:\spZ\to(-\infty, \infty]$ is
  the \textit{generalized Moreau enhancement (GME)}~\cite{LiGME} of a seed convex function $\Psi$:
\begin{equation}
  \label{eq:GME}
  \Psi_B(\cdot) \coloneqq \Psi(\cdot) - \min_{v\in\spZ}\left[
    \Psi(v) + \frac{1}{2}\norm{B(\cdot - v)}^2_{\sptildeZ}
  \right].
\end{equation}
(Note: The GME function $\Psi_B$ in \eqref{eq:cLiGME-w-general-fidelity} with $B$ being a zero matrix reproduces the seed convex function $\Psi$.)
\end{problem}
We formulate the proposed model \eqref{eq:cLiGME-w-general-fidelity} motivated by recent studies \cite{zhang2010,selesnick2017,LiGME,shabili2021,cLiGME,zhang2023} on the \textit{convexity-preserving nonconvex regularizers} (see "Related works" in the end of this section).
Indeed, the model \eqref{eq:cLiGME-w-general-fidelity} is an extension of the \textit{(constrained) Linearly involved Generalized Moreau Enhanced (LiGME) model} \cite{LiGME, cLiGME} in order to use a more flexible convex data fidelity term than the quadratic data fidelity case $f\coloneqq\frac{1}{2}\norm{y - \cdot}_{\spY}^2$ where (i) the cost function $J$ in \eqref{eq:cLiGME-w-general-fidelity} achieves the convexity by a strategical choice \cite{Chen2023} of a tuning matrix $B$ and (ii) inner-loop free algorithms were given with convergence guarantees to a global minimizer of the model \eqref{eq:cLiGME-w-general-fidelity}. 
For broader applications such as computed tomography \cite{bouman1996}, an extension of $f$ in the model in \cite{cLiGME} has also been studied in \cite{zhang2023} together with (i) the overall convexity condition of the extended model therein and (ii) an applicable difference-of-convex (DC) type algorithm, along the standard strategy in the DC programming \cite{tao1997, thi2018}, which requires inner loops.

The goal of this paper is to establish an inner-loop free algorithm for the model \eqref{eq:cLiGME-w-general-fidelity} by extending proximal splitting type algorithms in \cite{LiGME,cLiGME}.
More precisely, we address the following research questions according to the model \eqref{eq:cLiGME-w-general-fidelity} and its algorithm:
\begin{enumerate}[(Q1)]
  \item Under what conditions can we guarantee the existence of a minimizer of the model \eqref{eq:cLiGME-w-general-fidelity} ?
  \item How can we extend the existing proximal splitting type algorithms, developed specially for the quadratic data fidelity case \cite{LiGME,cLiGME}, to be applicable to the extended model \eqref{eq:cLiGME-w-general-fidelity}~?
\end{enumerate}

The remainder of this paper is organized as follows.
In Section \ref{sec:preliminaries}, as preliminaries, we introduce notation and selected tools in convex analysis and fixed point theory.
In Section~\ref{sec:main}, we first introduce a sufficient condition, called \textit{the overall convexity condition}, for the convexity of $J$ in \eqref{eq:cLiGME-w-general-fidelity} over $\spX$ (see Proposition \ref{prop:overall-convexity}).
Under the overall convexity condition, we present in Theorem \ref{thm:existence-solution} sufficient conditions for the existence of a minimizer using a classical theorem by Auslender \cite{Auslender1996}.
Under the overall convexity condition and the existence of a minimizer, we also propose a proximal splitting type algorithm with guaranteed convergence to a global minimizer of \eqref{eq:cLiGME-w-general-fidelity} (see Theorem \ref{thm:convergence-analysis}). 
The proposed algorithm requires Lipschitz continuity of $\nabla f$ over $\spY$.
However, for some applications (e.g., see Section \ref{sec:experiments}), Lipschitz continuity of $\nabla f$ can be assumed only over the constraint set $\opC^{-1}(\mathbf{C})$.
As a remedy of this issue, under a certain assumption, we present a way to reformulate such an optimization model into a model to which the proposed algorithm can be applied (see Proposition \ref{prop:convex-extension}).
In Section \ref{sec:experiments}, as an application of the proposed methods, we enhance nonconvexly a convex sparsity-promoting regularizer in a convex model \cite{banerjee2024} for simultaneous declipping and denoising.

\subsection*{Related works}
\label{sec:related-works}
{
The \textit{convexity-preserving nonconvex regularizers} was pioneered by Blake and Zisserman \cite{blake1987} and by Nikolova \cite{nikolova1998, nikolova1999}, and has been developed further.
The most of convexity-preserving nonconvex regularizers depend on the strong convexity of the least squares term $\frac{1}{2}\norm{y- A\cdot}_{\spY}^2$, i.e., the nonsingularity of $A^*A$.
As an exceptional example which is free from the nonsingularity of $A^*A$, Selesnick \cite{selesnick2017} proposed the \textit{generalized minimax concave (GMC) penalty} as a nonconvex enhancement of the $\ell_1$ norm\footnote{
  The $\ell_1$ norm $\norm{\cdot}_1$ is the largest convex minorant of the $\ell_0$ pseudo-norm in the vicinity of the zero vector.
  Thus, many models have used $\norm{\cdot}_1$ as a sparsity-promoting function. See, e.g., \cite{tibshirani1996,elad2010}.
} $\norm{\cdot}_1$ ($\Psi_B\circ \opL$ reproduces the GMC penalty with $\Psi = \norm{\cdot}_1$ and $\opL= \Id$). The GMC penalty is a nonseparable multidimensional extension of the minimax concave penalty \cite{zhang2010}.
To extend the idea of the GMC penalty for a general seed convex function $\Psi$ and a linear operator $\opL$, the LiGME regularizer $\Psi_B\circ \opL$ in \eqref{eq:GME} has been proposed~\cite{LiGME}.
See, e.g., \cite{feng2020,shabili2021,heng2025} for further advancements of the GMC model and \cite{kitahara2021,cLiGME,zhang2023,yata2024,kuroda2024,shoji2025} for those\footnote{
  Applications of LiGME type nonconvex enhancement are not limited to sparsity aware signal estimation. For example, it is applicable to discrete valued signal estimations \cite{shoji2025}.
} of the LiGME model.
}%
\section{Preliminaries}
\label{sec:preliminaries}
\subsection{Notation}
{
\thickmuskip=0\thickmuskip
\medmuskip=0\medmuskip
\thinmuskip=0.0\thinmuskip
\arraycolsep=0.5\arraycolsep
\label{sec:notation}

Symbols $\setN_0$, $\setR$, $\setNNR$ and $\setPR$ denote respectively all nonnegative integers, all real numbers, all nonnegative real numbers and all positive real numbers.
Let $\spH$ and $\spK$ be finite dimensional real Hilbert spaces.
A Hilbert space $\spH$ is equipped with an inner product $\ip{\cdot}{\cdot}_\spH$
and its induced norm $\norm{\cdot}_\spH$.
$\setLO{\spH}{\spK}$ denotes the set of all linear operators from 
$\spH$ to $\spK$.
For $L\in \setLO{\spH}{\spK}$, 
$\norm{L}_{\mathrm{op}}$ denotes the operator norm of $L$
(i.e., $\norm{L}_{\mathrm{op}}
:=\sup_{x\in\spH, \norm{x}_{\spH}\leq 1} \norm{Lx}_\spK$) and
$L^*\in\setLO{\spK}{\spH}$ the adjoint operator of $L$ 
(i.e., $(\forall x\in \spH)(\forall y \in\spK)\, 
\ip{Lx}{y}_\spK = \ip{x}{L^*y}_\spH$).
The identity operator is denoted by
$\Id$ and the zero operator from $\spH$ to $\spK$ by $\zeroMatrix_{\spH,\spK}$. In particular, we use the simplified notation $\zeroMatrix_{\spH}$ for the zero operator from $\spH$ to $\spH$.
We express the positive definiteness and the positive semidefiniteness of a self-adjoint operator $L\in \setLO{\spH}{\spH}$ 
as $L\succ \zeroMatrix_\spH$ and $L\succeq  \zeroMatrix_\spH$, respectively. 
Any $L\succ \zeroMatrix_{\spH}$ defines a new real Hilbert space $(\spH, \ip{\cdot}{\cdot}_{L}, \norm{\cdot}_L)$ where $\ip{x_1}{x_2}_{L}\coloneqq \ip{x_1}{L x_2}_{\spH}$ and $\norm{x}_{L}\coloneqq \sqrt{\ip{x}{L x}_{\spH}}$.
For a linear operator $L\in\setLO{\spH}{\spK}$ and sets $S_{\spH}\subset\spH$ and $S_{\spK}\subset\spK$, $L(S_{\spH})\coloneqq \{Lx \in\spK\mid x\in S_{\spH}\}$ is the image of $S_{\spH}$ under $L$, and $L^{-1}(S_{\spK})\coloneqq \{x\in\spH \mid Lx\in S_{\spK}\}$ is the preimage of $S_{\spK}$ under $L$.
A set $K$ is a cone if $(\forall x \in K)(\forall \alpha\in\setPR)\ \alpha x \in K$.
A set $C\subset \spH$ is convex if 
$(\forall x\in \spH)(\forall u\in\spH)(\forall \alpha\in [0,1])\ \alpha x +(1-\alpha) u \in C$.
A closed convex set $C\subset \spH$ is said to be \textit{simple} if its metric projection $P_C:\spH\to\spH: x\to\argmin_{u\in C} \norm{x-u}_{\spH}$ is available as a computable operator.
A function $f:\spH\to (-\infty,\infty]$ is said to be 
(i) proper if $\dom f:= \{x\in\spH \mid f(x)<\infty\}\neq \emptyset$, 
(ii) lower semicontinuous if $ \operatorname{lev}_{\leq\alpha}f \coloneqq \{ x\in \spH | f(x)\leq \alpha\}$ is closed for every $\alpha\in \setR$,
(iii) convex if $f(\alpha x + (1-\alpha) y)\leq \alpha f(x)+(1-\alpha) f(y) $ for every $x,y \in \spH,\, 0<\alpha< 1$.
The set of all proper lower semicontinuous convex functions defined on $\spH$ is denoted by $\Gamma_0 (\spH)$. $f\in\Gamma_0(\spH)$ is said to be coercive if $\lim_{\norm{x}_{\spH}\to\infty} f(x) = \infty$.
$f\in\Gamma_0(\spH)$ is said to be prox-friendly if 
$\prox_{\gamma f}:\spH\to\spH:x\mapsto \argmin_{v\in\spH}\left[
  f(v) + \frac{1}{2\gamma}\norm{v-x}^2_\spH
\right]$
is available as a computable operator for every $\gamma>0$.
For a nonempty closed convex set $C\subset\spH$, the proximity operator of the indicator function $\iota_C\in\Gamma_0(\spH)$ defined by 
$
  \iota_C(x) \coloneqq\begin{cases}
    0 & x\in C\\
    \infty & x\notin C
  \end{cases}
$
coincides with the metric projection operator onto $C$, i.e., $\prox_{\gamma\iota_{C}} = P_C \ (\forall \gamma > 0)$.
}%
\subsection{Selected tools in fixed point theory}
An operator $T:\spH\to\spH$ is said to be \textit{nonexpansive} if
\begin{equation}
  (\forall x\in\spH)(\forall y\in\spH)\ 
  \norm{Tx-Ty}_{\spH} \leq \norm{x-y}_{\spH}.
\end{equation}
In particular, $T$ is \textit{$\alpha$-averaged nonexpansive} with $\alpha\in (0,1)$ if there exists a nonexpansive operator $\widehat{T}:\spH\to\spH$ such that $T= (1-\alpha)\Id + \alpha \widehat{T}$.

A fixed point of a nonexpansive operator can be approximated successively via the Krasnosel'ski\u{\i}-Mann iteration.
\begin{fact}[Krasnosel'ski\u{\i}-Mann iteration (See, e.g., {\cite{groetsch1972}})]
  \label{fact:KM-itr}
  Let $T:\spH\to\spH$ be a nonexpansive operator such that $\fix(T) := \{x \in \spH\mid x = T(x)\}\neq \emptyset$.
    For any initial point $x_0\in\spH$,
    set $(x_k)_{k\in\setN_0}$ by 
    \begin{equation}
      (\forall k\in \setN_0)\  x_{k+1} = [(1-\lambda_k)\Id + \lambda_k T](x_k)
    \end{equation}
    with a sequence $(\lambda_k)_{k\in\setN_0}$ 
    satisfying $(\forall k\in \setN_0)\ \lambda_k\in [0,1]$ 
    and
    $\sum_{k\in\setN_0}\lambda_k(1- \lambda_k)=+\infty$.
    Then $(x_k)_{k\in\setN_0}$ converges to a point in $\fix(T)$.
    In particular, if $T$ is $\alpha$-averaged with $\alpha\in(0,1)$,
    the sequence generated by 
    $
      (\forall n \in \setN_0)\ x_{k+1} =T(x_k)
    $
    converges to a point in $\fix(T)$.
\end{fact}

\section{Properties of Problem \ref{prob:cLiGME-w-general-fidelity} and iterative algorithm for wider applications}
\label{sec:main}
\subsection{Overall convexity and existence of minimizer}
We start with a sufficient condition for the convexity of the proposed cost function $J$ in \eqref{eq:cLiGME-w-general-fidelity}.

\begin{proposition}
  \label{prop:overall-convexity}
  In Problem \ref{prob:cLiGME-w-general-fidelity}, set
\begin{align}
  \mathfrak{d} :\spX\to\setR: x\mapsto f\circ A(x) - \frac{\mu}{2}\norm{B\opL x}^2_{\sptildeZ}.
\end{align}
Then the relation
    $(C_1)$: $ \mathfrak{d}\in\Gamma_0(\spX)\implies
    (C_2)$: $J\in\Gamma_0(\spX)$ holds.
    In particular, if $f$ is twice continuously differentiable over $\spY$, then the condition $(C_1)$ is equivalent to
    \begin{equation}
      \label{eq:C-1dash}
      (C_1'): (\forall x\in\spX)\ A^*\circ \nabla^2 f(Ax)\circ A- \mu \opL^*B^*B\opL\succeq \zeroMatrix_{\spX}.
    \end{equation}
    We call $(C_1)$ the overall convexity condition.
\end{proposition}
\begin{remark}[Comparison with existing conditions]
  For wider applicable conditions than $(C_1)$ and $(C_1')$, see \cite[Theorem~1 and Corollary 1]{zhang2023}.
  In a quadratic case $f\coloneqq\frac{1}{2}\norm{\cdot -y}^2_{\spY}$,
  the condition $(C_1')$ reproduces the condition in \cite[Proposition 1]{LiGME} which is a generalization of a special case in \cite[Theorem 1]{selesnick2017} applicable for $(\Psi,\opL)\coloneqq (\norm{\cdot}_1,\Id)$.
\end{remark}
\begin{remark}[An algebraic design of $B$ enjoying $(C_1')$]
  \label{remark:choice-B}
 Given $\Lambda\succeq\zeroMatrix_{\spY}$ such that 
  \begin{equation}
    \label{eq:lambda-lower-bound}
    (\forall x\in\spX) \ \nabla^2f(Ax)-\Lambda\succeq\zeroMatrix_{\spY},
  \end{equation}
  design $B$ enjoying
  \begin{equation}
    \label{eq:lambda-oc-cond}
    A^*\Lambda A-\mu\opL^*B^*B\opL\succeq\zeroMatrix_{\spX}.
  \end{equation}
  GME matrix $B$ enjoying \eqref{eq:lambda-oc-cond} can be designed
  via LDU decomposition of $\opL$ \cite[Theorem 1]{Chen2023}.
  Then $B$ achieves the condition $(C_1')$.
  In particular, if $\spY\coloneqq \setR^m$ and $f\in\Gamma_0(\setR^m)$ is a separable sum $f(u)\coloneqq \sum_{i=1}^m f_i([u]_i)$ of twice continuously differentiable convex functions $f_i\in\Gamma_0(\setR)$ on real line, then a diagonal matrix $\Lambda\in\setR^{m\times m}$ with $i$-th diagonal entry
  \begin{equation}
    \label{eq:choice-Lambda}
    [\Lambda]_{i,i} = \inf_{r\in\setR} f''_i(r)\geq 0
  \end{equation}
  satisfies the condition \eqref{eq:lambda-lower-bound}.
\end{remark}
Essentially based on \cite[Theorem 2.4]{Auslender1996},
we derived the following sufficient conditions to guarantee the existence of a minimizer of the model \eqref{eq:cLiGME-w-general-fidelity} under the overall convexity condition $(C_1)$.
\begin{theorem}[Existence of a minimizer]
  \label{thm:existence-solution}
  Consider Problem \ref{prob:cLiGME-w-general-fidelity}.
  Assume that the condition $(C_1)$ in Proposition \ref{prop:overall-convexity} is achieved.
  Then the model \eqref{eq:cLiGME-w-general-fidelity} has a minimizer if one of the following holds.
  \begin{enumerate}[\upshape (i)]
  \item $f$ is coercive, and the constraint set $\opC^{-1}(\mathbf{C})$ is asymptotically multipolyhedral \cite[Definition 2.3]{Auslender1996}, i.e., the constraint set $\opC^{-1}(\mathbf{C})$ can be decomposed as $\opC^{-1}(\mathbf{C})= S+K$ with a compact set $S\subset\spX$ and a polyhedral cone $K\subset\spX$, where the addition of two sets is understood as the Minkowski sum.
  \item The constraint set $\opC^{-1}(\mathbf{C})$ is bounded.
    \item $f$ is coercive and $\nullsp A\cap \nullsp \opL = \{0_{\spX}\}$.
  \end{enumerate} 
\end{theorem}
For a condition for $\opC^{-1}(\mathbf{C})$ to admit such a decomposition $\opC^{-1}(\mathbf{C})= S+K$ in (i), see, e.g., \cite{Goberna2010}.
The asymptotically multipolyhedral convex set in the condition (i) of Theorem~\ref{thm:existence-solution} can cover a wide range of convex sets as follows.
\begin{example}[Asymptotically multipolyhedral sets]
  \quad
  \begin{enumerate}[(a)]
    \item (Polyhedral set). Any polyhedral set can be decomposed as the sum of a polytope (which is compact) and a polyhedral convex cone \cite{Goberna2010}. Therefore, any polyhedral convex set is asymptotically multipolyhedral.
    \item (Entire space). The entire space $\spX$ is a typical example of polyhedral sets and thus asymptotically multipolyhedral as well. 
    For the LiGME model, i.e., the model \eqref{eq:cLiGME-w-general-fidelity} with $f=\frac{1}{2}\norm{y- \cdot}_{\spY}^2$ and $\opC^{-1}(\mathbf{C})=\spX$, the existence of its minimizer is  guaranteed by the condition (i) in Theorem 1, while such an existence is assumed implicitly in \cite{LiGME}.
    \item (Linearly involved compact set). Let $\mathbf{C}$ be a compact convex set.
    Then $\opC^{-1}(\mathbf{C})$ can be decomposed as 
    $\opC^{-1}(\mathbf{C}) =  \opC^\dag (\mathbf{C}) + \ker \opC$, where $\opC^\dag$ is the Moore-Penrose pseudo inverse of $\opC$, and hence $\opC^\dag (\mathbf{C})$ is compact.
  \end{enumerate}
\end{example}

\subsection{Proximal splitting type algorithm for Problem \ref{prob:cLiGME-w-general-fidelity} with guaranteed convergence to a global minimizer}
\label{sec:algorithm}
We propose an iterative  algorithm for finding a global minimizer of the model \eqref{eq:cLiGME-w-general-fidelity} under the following assumption.

\begin{assumption}
  \label{assumption:alg}
  In Problem \ref*{prob:cLiGME-w-general-fidelity}, assume the following.

  \begin{enumerate}[(i)]
    \item $\nabla \mathfrak{d} = A^*\circ \nabla f\circ A - \mu\opL^*B^*B\opL$ is $\beta$-Lipschitz continuous over $\spX$ for some $\beta\in\setPR$.
    \item The overall convexity condition $(C_1)$ in Proposition \ref{prop:overall-convexity} is satisfied (see Remark \ref{remark:choice-B} for choice of $B$ enjoying $(C_1)$).
    \item A minimizer of \eqref{eq:cLiGME-w-general-fidelity} exists (see Theorem \ref{thm:existence-solution} for sufficient conditions).
    \item $\Psi\circ\opL$, and $\iota_{\mathbf{C}}\circ \opC$ in \eqref{eq:cLiGME-w-general-fidelity} satisfy\footnote{
      By \cite[Proposition 6.19, Corollary 16.50 and Corollary 16.53]{CAaMOTiH},
  Assumption \ref{assumption:alg}(iv) is satisfied if 
  $
    \emptyset\neq \ri(\dom(\Psi\circ \opL))\cap \opC^{-1}(\mathbf{C}) \mbox{ and } \emptyset\neq \ri (\mathbf{C}) \cap \ran \opC.
  $
  (For a convex set $C$, $\ri C$ denotes the relative interior of $C$. See, e.g., \cite[Definition 6.9]{CAaMOTiH}.)
    }
    \begin{equation}
      \partial(\mu\Psi\circ\opL +  \iota_{\mathbf{C}}\circ \opC)
      = \mu\opL^*\circ(\partial \Psi)\circ\opL + \opC^*\circ(\partial\iota_{\mathbf{C}})\circ \opC,
    \end{equation}
    where the subdifferential $\partial \Psi:\spZ\to 2^\spZ$ of $\Psi\in\Gamma_0(\spZ)$ is defined as 
    {
      \thickmuskip=0\thickmuskip
      \medmuskip=0\medmuskip
      \thinmuskip=0.0\thinmuskip
      \arraycolsep=0.5\arraycolsep
    \begin{equation}
      \partial\Psi(z) \coloneqq\{u\in\spZ\mid (\forall w \in\spZ)\ \ip{w-z}{u}_{\spZ} + \Psi(z)\leq \Psi(w)\}.
    \end{equation}
    }
  \end{enumerate}
\end{assumption}
The following theorem shows that the set of all minimizers of the model \eqref{eq:cLiGME-w-general-fidelity} can be expressed in terms of the fixed point set of an averaged nonexpansive operator $T$ in Theorem \ref{thm:convergence-analysis}, and thus a global minimizer can be approximated iteratively by the Krasnosel'ski\u{\i}-Mann iteration of $T$.
\begin{theorem}[Fixed point translation of the solution set of \eqref{eq:cLiGME-w-general-fidelity}]
  \label{thm:convergence-analysis}
  Consider Problem \ref{prob:cLiGME-w-general-fidelity} under Assumption \ref{assumption:alg}. Let $\setS$ be the set of all global minimizers of the model \eqref{eq:cLiGME-w-general-fidelity}.
  Set a product space $\spH:= \spX\times \spZ\times\spZ\times \spfrakZ$ and define an operator:
  \begin{equation}
    T:\spH\to\spH: (x,v,w,z)\mapsto (\xi,\zeta,\eta,\varsigma)
  \end{equation}
  with $(\sigma,\tau)\in \setPR\times\setPR$ by
  \begin{align}
    \xi &:= \left(\Id - \frac{1}{\sigma} \nabla \mathfrak{d}\right) (x) - \frac{\mu}{\sigma}\opL^*B^*B v -\frac{\mu}{\sigma}\opL^*w - \frac{\mu}{\sigma}\opC^*z\\
    \zeta &:= \prox_{\frac{\mu}{\tau}\Psi}\left[
      \frac{2\mu}{\tau}B^*B\opL\xi
      -\frac{\mu}{\tau} B^*B\opL x + \left(\Id -\frac{\mu}{\tau}B^*B \right)(v)
    \right]\\
    \eta&:= (\Id - \prox_{\Psi})(2\opL \xi -\opL x +w)\\
    \varsigma &:= (\Id - P_{\mathbf{C}})(2\opC\xi - \opC x +z).
  \end{align}
  Then the following hold.
  \begin{enumerate}[\upshape (a)]
    \item The solution set $\mathcal{S}$ of \eqref{eq:cLiGME-w-general-fidelity} can be expressed as 
    \begin{align}
      \label{eq:fixed-point-encode}
      \mathcal{S}=\Xi(\fix T):=\{
        \Xi(h) \in\spX\mid 
    h = T(h)
    \},
    \end{align}
    where $\Xi:\spH\to\spX:(x,v,w,z)\mapsto x$.
    \item Choose $(\sigma, \tau)\in(0,\infty)\times (\frac{1}{2\rho},\infty)$ satisfying
    \begin{align}
      \label{eq:cond-sigma-tau}
      \sigma > \mu \norm{\opL^*\opL + \opC^*\opC}_{\mathrm{op}}+\frac{2\rho \mu^2\norm{B^*B\opL}^2_{\mathrm{op}} + \tau}{2\rho\tau - 1},
    \end{align}
    where $\rho \coloneqq\frac{1}{\max\left\{  \beta, \mu\norm{B}_{\mathrm{op}}^2\right\}}>0$.
    Then 
    \begin{equation}
      \label{eq:def-P}
      \opP \coloneqq \begin{bmatrix}
        \sigma \Id & -\mu\opL^*B^*B & -\mu \opL^* & -\mu\opC^*\\
        -\mu B^* B \opL & \tau \Id & \zeroMatrix_{\spZ, \spZ} & \zeroMatrix_{\spfrakZ,\spZ}\\
        -\mu \opL & \zeroMatrix_{\spZ,\spZ} & \mu \Id & \zeroMatrix_{\spfrakZ,\spZ}\\
        -\mu\opC & \zeroMatrix_{\spZ,\spfrakZ} & \zeroMatrix_{\spZ,\spfrakZ} & \mu\Id
      \end{bmatrix}\succ \zeroMatrix_{\spH}.
    \end{equation}
    Furthermore,
    $T$ is $\frac{2}{4-\theta}$-averaged nonexpansive over $(\spH, \ip{\cdot}{\cdot}_{\opP}, \norm{\cdot}_{\opP})$ with
    \begin{align}
      \theta:= \frac{ \sigma  + \tau- \mu\norm{\opL^*\opL + \opC^*\opC }_{\mathrm{op}}}{\rho(\sigma \tau - \tau\mu\norm{\opL^*\opL + \opC^*\opC }_{\mathrm{op}} - \mu^2 \norm{B^* B \opL}_{\mathrm{op}}^2)} \in (0, 2).
    \end{align}
    \item For any initial point $h_0\in\spH$, the sequence $(h_k)_{k\in\setN_0} \subset \spH$ generated by the Krasnosel'ski\u{\i}-Mann iteration:
    \begin{equation}
      \label{eq:K-M}
      (\forall k \in\setN_0)\ h_{k+1} \coloneqq T(h_k)
    \end{equation} 
    converges to a point in $\fix T$ (see Fact \ref{fact:KM-itr}), which implies that the sequence $(\spX\ni)x_k \coloneqq \Xi (h_k) \ ( k \in\setN_0)$ converges to a point in the solution set $\setS$.
    \end{enumerate}
\end{theorem}

\begin{remark}[Comparision with existing algorithms]
  The proposed algorithm can be seen as an extension of \cite[Algorithm 1]{LiGME} and \cite[Algorithm 1]{cLiGME} which are proposed for the quadratic data fidelity case, i.e., $f\coloneqq \frac{1}{2}\norm{y-\cdot}^2_{\spY}$.
  Moreover, the proposed algorithm is applicable to general cases where the conditions $\dom\Psi = \spZ$ and $\Psi\circ(-\Id) = \Psi$, imposed in \cite[Algorithm 1]{LiGME} and \cite[Algorithm 1]{cLiGME}, are no longer satisfied.
\end{remark}

\subsection{Remedy for lack of Lipschitz continuous gradient of data fidelity functions}
\label{sec:extension}
For some applications, Lipschitz continuity assumption in Assumption \ref{assumption:alg}(i) does not hold (see, e.g., the model \eqref{eq:numerical-proposed} in Section \ref{sec:experiments}).
For such applications, we consider the model~\eqref{eq:cLiGME-w-general-fidelity} under Assumption \ref{assumption:alg}(ii)-(iv) and the following alternative assumption against Assumption \ref{assumption:alg}(i).
Hereafter, we denote the $i$-th component of a vector $u\in\setR^m$ by $[u]_i\in\setR$.

\begin{assumption}[Alternative assumption against Assumption~\ref{assumption:alg}(i)]
  \label{assumption:non-Lipschitz}
  In Problem \ref*{prob:cLiGME-w-general-fidelity}, let $\spY\coloneqq \setR^m$ and a closed convex set $\varPi \supset A(\opC^{-1}(\mathbf{C}))$
  be decomposable as
  $
    \varPi \coloneqq \bigtimes_{i=1}^m \varPi_i (\subset \spY)
  $
  with a closed convex set $\varPi_i\subset\setR$.
  Assume that $f \in\Gamma_0(\setR^m)$ is a separable sum:
  \begin{equation}
    \label{eq:separable-assumption}
    (\forall u\in \mathbb{R}^{m})\ 
    f(u) = \sum_{i=1}^m f_i ([u]_i)
  \end{equation}
  of $f_i\in\Gamma_0(\setR)\ (1\leq i \leq m)$, where $f_i\in\Gamma_0(\setR)\ (1\leq i \leq m)$ are twice continuously differentiable over $\varPi_i$ and
  $f_i'\ (1\leq i \leq m)$ are Lipschitz continuous over $\varPi_i$ (not necessarily Lipschitz continuous over $\setR$).
\end{assumption}

Due to the lack of Lipschitz continuity of $\nabla f$ over $\spY$, 
Theorem \ref{thm:convergence-analysis}(c) under Assumption \ref{assumption:alg}(ii)-(iv) and Assumption \ref{assumption:non-Lipschitz} no longer guarantees that the proposed algorithm produces a convergent sequence to a global minimizer of the model \eqref{eq:cLiGME-w-general-fidelity}.
To circumvent this issue, we introduce an alternative function $\widetilde{f}\in\Gamma_0(\spY)$, of $f$, enjoying desired properties as a data fidelity function applicable to the proposed algorithm.

\begin{proposition}[Construction of alternative data fidelity function $\widetilde{f}$ of $f$]
  \label{prop:convex-extension}
  Consider an instance of the model \eqref{eq:cLiGME-w-general-fidelity} in Problem \ref{prob:cLiGME-w-general-fidelity}
  satisfying Assumption \ref{assumption:alg}(ii)-(iv) and Assumption \ref{assumption:non-Lipschitz}.
Then we have the following.
\begin{enumerate}[\upshape (a)]
  \item
  Define $\widetilde{f}:\setR^m\to\setR$ by $\widetilde{f}(u)\coloneqq \sum_{i=1}^m\widetilde{f}_i([u]_i)$ 
  with univariate functions $ (1\leq i \leq m) \ \widetilde{f_i}:\setR\to\setR:$ 
  {
  \begin{equation}
    \label{eq:def-gtilde-i}
    r\mapsto \begin{cases}
      \frac{f_i''(c_r)}{2} (r - c_r)^2 +  f_i'(c_r)
      (r-c_r) + f_i(c_r) & r\notin\varPi_i\\
      f_i(r) & r \in \varPi_i,
    \end{cases}
  \end{equation}
  }%
  where $c_r\coloneqq P_{\varPi_i} (r)\in\setR$.
  Then, for an optimization model:
  \begin{equation}
    \label{eq:extended-model}
    \minimize_{\opC x \in\mathbf{C}} \widetilde{f}\circ A(x) + \mu \Psi_B\circ\opL (x),
  \end{equation}
  we have
  {
    \thickmuskip=0\thickmuskip
    \medmuskip=0\medmuskip
    \thinmuskip=0.0\thinmuskip
    \arraycolsep=0.5\arraycolsep
  \begin{equation}
    \label{eq:argmin-equivalence}
    \argmin_{\substack{\opC x \in\mathbf{C}}} f\circ A(x) + \mu \Psi_B\circ\opL (x)
    = \argmin_{\opC x \in\mathbf{C}} \widetilde{f}\circ A(x) + \mu \Psi_B\circ\opL (x).
  \end{equation}
  }%
  Moreover, the model \eqref{eq:extended-model} is also an instance of Problem \ref{prob:cLiGME-w-general-fidelity} and enjoys the following conditions:
  \begin{align*}
    &(\forall x\in \spX) \ A^*\circ \nabla^2\widetilde{f}(Ax)\circ A - \mu\opL^*B^*B\opL\succeq\zeroMatrix_\spX,\\
    &\nabla \left(\widetilde{f}\circ A - \frac{\mu}{2}\norm{B\opL \cdot}^2_{\sptildeZ}\right) \mbox{ is Lipschitz continuous over } \spY.
  \end{align*}
  \item 
  Generate the sequence $(h_k)_{k\in\setN_0}$ by the proposed algorithm~\eqref{eq:K-M} with a replacement of $\mathfrak{d}$ with 
  $\widetilde{\mathfrak{d}}\coloneqq \widetilde{f}\circ A - \frac{\mu}{2}\norm{B\opL \cdot}^2_{\sptildeZ}$.
  Then, by Theorem \ref{thm:convergence-analysis}, the sequence $x_k\coloneqq\Xi(h_k)\ (k\in\setN_0)$ converges to a global minimizer $\bar{x}\in\opC^{-1}(\mathbf{C})$ of the model~\eqref{eq:extended-model}.
  Moreover, from \eqref{eq:argmin-equivalence}, such a minimizer $\bar{x}\in\opC^{-1}(\mathbf{C})$ is also a minimizer of the model~\eqref{eq:cLiGME-w-general-fidelity}.
\end{enumerate}
\end{proposition}

\section{Application to simultaneous declipping and denoising}
\label{sec:experiments}
\subsection{Formulation via Problem \ref{prob:cLiGME-w-general-fidelity}}
As an application of the proposed model \eqref{eq:cLiGME-w-general-fidelity}, we consider a problem for simultaneous declipping and denoising.
The task is to estimate the target signal $\xstar\in\setR^m$ from
\begin{equation}
  \label{eq:clip-observation}
  y = \operatorname{clip}_{\vartheta}(\xstar+\varepsilon) \in\setR^m
\end{equation} 
with a priori knowledge that (i) $\xstar \in \mathbf{C}$ with a nonempty closed convex set $\mathbf{C}$ and (ii) $\opL \xstar$ is sparse with a certain linear operator $\opL$, where $\operatorname{clip}_{\vartheta}:\setR^m \to \setR^m$ is defined with $\vartheta\in\setPR$ as an entrywise operator:
\begin{equation}
  [\operatorname{clip}_{\vartheta} (u)]_i \coloneqq \begin{cases}
    [u]_i & |[u]_i| < \vartheta\\
    \vartheta \cdot\operatorname{sign}([u]_i) & |[u]_i| \geq\vartheta
  \end{cases}
\end{equation}
for every $u\in\setR^m$ and $i\in\{1,2,\ldots,m\}$ 
and the noise values $[\varepsilon]_i$ follows Gaussian distribution with zero mean and known variance $s^2\ (s>0)$. 
Recently, \cite{banerjee2024} has proposed a data fidelity function $f:\setR^m\to\setR: u\mapsto \sum_{i=1}^m f_i([u]_i)$ with
\begin{equation}
  f_i([x]_i)\coloneqq
  \begin{cases}
    \frac{1}{2}\left(\frac{[y]_i - [x]_i}{s}\right)^2 & -\vartheta < [y]_i < \vartheta\\
    -\log \left(\int_{\vartheta - [x]_i}^{\infty}\exp\left(
      - \frac{t^2}{2s^2}
    \right) \mathrm{d}t \right) &[y]_i = \vartheta\\
    -\log \left(\int_{-\infty}^{-\vartheta- [x]_i}\exp\left(
      - \frac{t^2}{2s^2}
    \right) \mathrm{d}t \right) &[y]_i = -\vartheta
  \end{cases}
\end{equation}
for the observation model \eqref{eq:clip-observation}, and formulated\footnote{The original model in \cite{banerjee2024} employed $f\circ A$ with a linear operator $A$ as a data fidelity term and was given without constraints.
To clarify the effectiveness of nonconvex enhancement, we consider the simplest case where $A\coloneqq\Id$ although the proposed model and algorithm can also be applied to a general case $A\neq\Id$.} an optimization model:
\begin{equation}
  \label{eq:numerical-conventional}
  \minimize_{x\in\mathbf{C}} f(x) + \mu \norm{\cdot}_{1}\circ \opL (x),
\end{equation}
where $\norm{\cdot}_{1}\in\Gamma_0(\setR^{m})$ is the $\ell_1$ norm\footnote{$\norm{\cdot}_{1}$ is prox-friendly. See, e.g., \cite[Example 24.22]{CAaMOTiH}.}.
To enhance nonconvexly $\norm{\cdot}_1$ in \eqref{eq:numerical-conventional}, we propose 
\begin{equation}
  \label{eq:numerical-proposed}
  \minimize_{x\in\mathbf{C}} f(x) + \mu (\norm{\cdot}_{1})_B\circ \opL (x)
\end{equation}
as a special instance of the model~\eqref{eq:cLiGME-w-general-fidelity}.
By employing a GME matrix $B$ enjoying $(C_1)$ in Proposition \ref{prop:overall-convexity} and a compact set as $\mathbf{C}$, the model \eqref{eq:numerical-proposed} enjoys Assumption \ref{assumption:alg}(ii-iv) and Assumption \ref{assumption:non-Lipschitz}.
Therefore, we can approximate iteratively a minimizer of \eqref{eq:numerical-proposed} by applying the proposed algorithm \eqref{eq:K-M} to an alternative optimization model~\eqref{eq:extended-model} in Proposition~\ref{prop:convex-extension}.

\subsection{Numerical experiments}
Following \cite{banerjee2024}, we conducted numerical experiments on estimation of $\xstar \in \setR^m \ (m\coloneqq256)$ from its noisy observation $y$ in~\eqref{eq:clip-observation}, where $\xstar$ was given by the inverse Discrete Cosine Transform (DCT) of a randomly chosen sparse coefficient vector. The target signal $\xstar$ was normalized to satisfy $\norm{\xstar}_{\infty}=0.8$.
We compared the conventional model \eqref{eq:numerical-conventional} and the proposed model \eqref{eq:numerical-proposed} for every $(\vartheta, s)\in\{0.4,0.6\}\times \{s_{5},s_{10},s_{15}\}$, where $s_5, s_{10}$ and $s_{15}$ are standard deviations of Gaussian noise achieving respectively $5$dB, $10$dB and $15$dB of SNR: $20\log_{10}\frac{\norm{x^\star}_{\mathbb{R}^m}}{\mathbb{E}[\norm{\varepsilon}_{\mathbb{R}^m}]}$.
For both models \eqref{eq:numerical-conventional} and \eqref{eq:numerical-proposed}, we employed $\mathbf{C}\coloneqq [-10,10]^m$ and $\opL\coloneqq \opL_{\mathrm{DCT}}$ (DCT matrix \cite{rao2007}), and a simple GME matrix\footnote{
  $\Lambda\in\setR^{m\times m}$ was given by \eqref{eq:choice-Lambda} with 
$
  \inf_{r\in\setR} f_i''(r) = \begin{cases}
    \frac{1}{s^2} & -\vartheta < [y]_i <\vartheta\\
    0 & |[y_i]|=\vartheta.
  \end{cases}
  $
} $B\coloneqq \sqrt{\frac{0.99}{\mu}} \sqrt{\Lambda} \opL_{\mathrm{DCT}}^{-1}$ for achieving the condition $(C_1)$ in Proposition \ref{prop:overall-convexity} for the model \eqref{eq:numerical-proposed}.
For minimization of \eqref{eq:numerical-conventional} and \eqref{eq:numerical-proposed}, we introduced an alternative data fidelity function $\widetilde{f}$ as in Proposition \ref{prop:convex-extension} with $\varPi\coloneqq\mathbf{C}$ and then applied\footnote{
  The model \eqref{eq:numerical-proposed} reproduces the model \eqref{eq:numerical-conventional} by setting a zero matrix to $B$. Thus, we can apply the proposed algorithm \eqref{eq:K-M} to the model \eqref{eq:numerical-conventional}.
}
the proposed algorithm\footnote{We used $\tau\coloneqq\frac{5}{2\rho}$ and $\sigma$ given by $1.001\times$(the value of RHS in \eqref{eq:cond-sigma-tau}).} \eqref{eq:K-M}. 
We stopped the proposed algorithm \eqref{eq:K-M} after the residual achieves $\norm{h_k-h_{k-1}}_{\spH}<10^{-4}$.

\begin{figure}[t]
  \centering
  \begin{minipage}[b]{0.46\linewidth}
    \centering
    \includegraphics[keepaspectratio, scale=0.11]
    {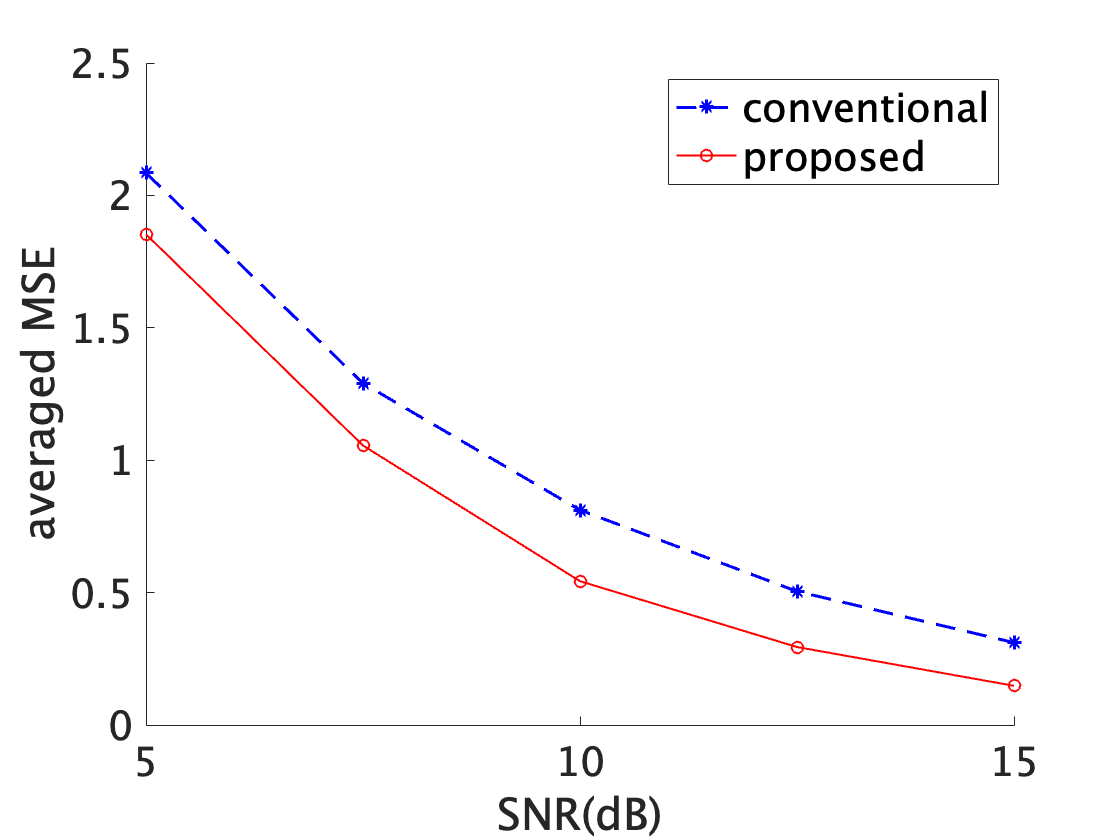}
    \subcaption{$\vartheta=0.4$}
  \end{minipage}
  \begin{minipage}[b]{0.46\linewidth}
    \centering
    \includegraphics[keepaspectratio,scale=0.11]
    {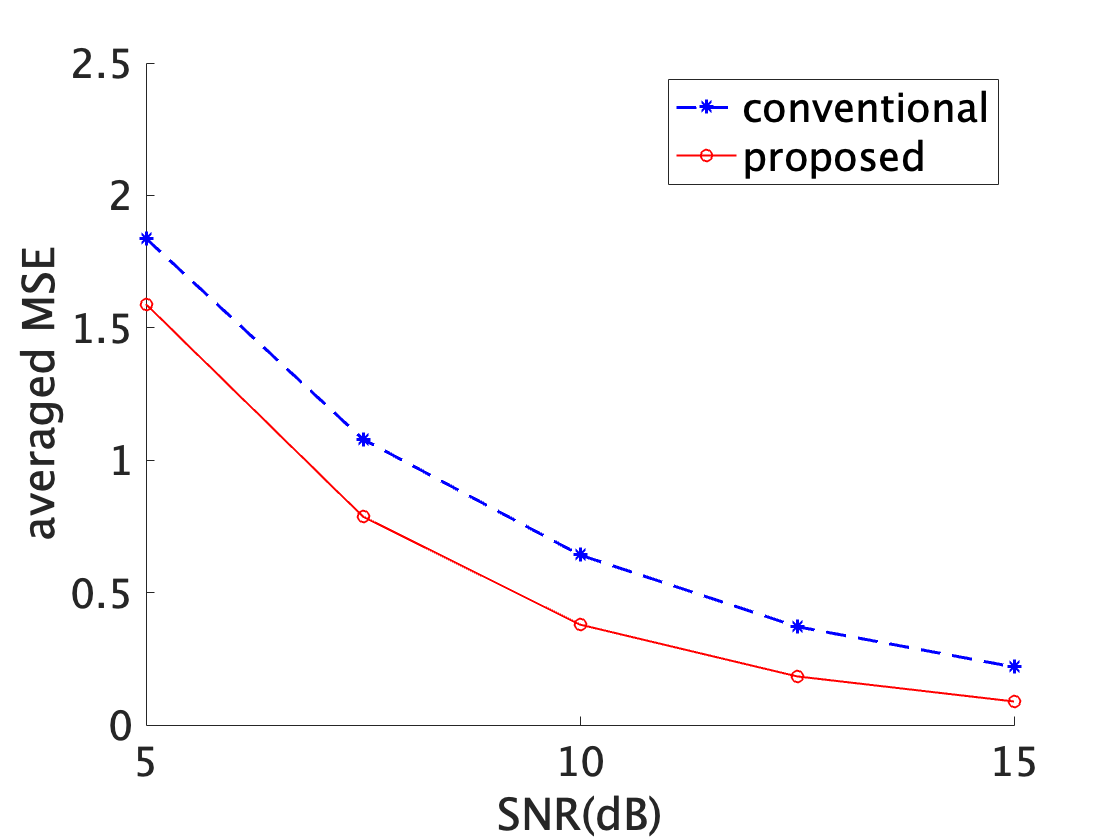}
    \subcaption{$\vartheta=0.6$}
  \end{minipage}
  \caption{
    SNR versus Averaged MSE
  }
  \label{fig:denoise-example}
  \vspace{-1em}
\end{figure}
Figure \ref{fig:denoise-example} shows the average of MSE: $\norm{x^\star - \bar{x}}_2^2$ of estimates $\bar{x}$ by the existing model \eqref{eq:numerical-conventional} and by the proposed model~\eqref{eq:numerical-proposed} over 100 realizations of Gaussian noise $\varepsilon$.
For each $(\vartheta, s^2)$, we choose the best regularization parameter $\mu$ from $\{j\in\setN_0\mid 1\leq j \leq 100\}$.
From this figure, we see that the proposed model~\eqref{eq:numerical-proposed} outperforms the model \eqref{eq:numerical-conventional} in all cases.

\section{Conclusion}
In this paper, we proposed a nonconvexly regularized convex model with a smooth data fidelity and the LiGME regularizer. 
Under the overall convexity condition, we propose sufficient conditions for existence of a minimizer for the proposed model. We also propose a proximal splitting type algorithm for finding a global minimizer of the proposed model.
Numerical experiments demonstrate the effectiveness of the proposed model and algorithm in a scenario of simultaneous declipping and denoising.

\noindent
\textbf{Acknowledgment}:
We thank Mr. Shih-Yang Lin in our laboratory for kind help on numerical experiments.

\bibliographystyle{IEEEtran}
\bibliography{references}

\begin{thebibliography}{10}
\providecommand{\url}[1]{#1}
\csname url@samestyle\endcsname
\providecommand{\newblock}{\relax}
\providecommand{\bibinfo}[2]{#2}
\providecommand{\BIBentrySTDinterwordspacing}{\spaceskip=0pt\relax}
\providecommand{\BIBentryALTinterwordstretchfactor}{4}
\providecommand{\BIBentryALTinterwordspacing}{\spaceskip=\fontdimen2\font plus
\BIBentryALTinterwordstretchfactor\fontdimen3\font minus
  \fontdimen4\font\relax}
\providecommand{\BIBforeignlanguage}[2]{{%
\expandafter\ifx\csname l@#1\endcsname\relax
\typeout{** WARNING: IEEEtran.bst: No hyphenation pattern has been}%
\typeout{** loaded for the language `#1'. Using the pattern for}%
\typeout{** the default language instead.}%
\else
\language=\csname l@#1\endcsname
\fi
#2}}
\providecommand{\BIBdecl}{\relax}
\BIBdecl

\bibitem{elad2010}
M.~Elad, \emph{Sparse and Redundant Representations}.\hskip 1em plus 0.5em
  minus 0.4em\relax Springer, 2010.

\bibitem{theodoridis2024}
S.~Theodoridis, \emph{{Machine Learning: From the Classics to Deep Networks,
  Transformers, and Diffusion Models}}, 3rd~ed.\hskip 1em plus 0.5em minus
  0.4em\relax Academic Press, 2024.

\bibitem{LiGME}
J.~Abe, M.~Yamagishi, and I.~Yamada, ``Linearly involved generalized {Moreau}
  enhanced models and their proximal splitting algorithm under overall
  convexity condition,'' \emph{Inverse Problems}, vol.~36, no.~3, 2020.

\bibitem{zhang2010}
C.-H. Zhang, ``{Nearly unbiased variable selection under minimax concave
  penalty},'' \emph{Annals of Statistics}, vol.~38, no.~2, 2010.

\bibitem{selesnick2017}
I.~Selesnick, ``{Sparse Regularization via Convex Analysis},'' \emph{IEEE
  Transactions on Signal Processing}, vol.~65, no.~17, 2017.

\bibitem{shabili2021}
A.~H. Al-Shabili, Y.~Feng, and I.~Selesnick, ``{Sharpening Sparse Regularizers
  via Smoothing},'' \emph{IEEE Open Journal of Signal Processing}, vol.~2,
  2021.

\bibitem{cLiGME}
W.~Yata, M.~Yamagishi, and I.~Yamada, ``{A constrained LiGME model and its
  proximal splitting algorithm under overall convexity condition},''
  \emph{Journal of Applied and Numerical Optimization}, vol.~4, no.~2, 2022.

\bibitem{zhang2023}
Y.~Zhang and I.~Yamada, ``{A Unified Framework for Solving a General Class of
  Nonconvexly Regularized Convex Models},'' \emph{IEEE Transactions on Signal
  Processing}, vol.~23, 2023.

\bibitem{Chen2023}
Y.~Chen, M.~Yamagishi, and I.~Yamada, ``{A Unified Design of Generalized Moreau
  Enhancement Matrix for Sparsity Aware LiGME Models},'' \emph{IEICE
  Transactions on Fundamentals of Electronics, Communications and Computer
  Sciences}, vol. E106.A, no.~8, 2023.

\bibitem{bouman1996}
C.~Bouman and K.~Sauer, ``A unified approach to statistical tomography using
  coordinate descent optimization,'' \emph{IEEE Transactions on Image
  Processing}, vol.~5, 1996.

\bibitem{tao1997}
P.~D. Tao and L.~T.~H. An, ``Convex analysis approach to dc programming:
  theory, algorithms and applications,'' \emph{Acta mathematica vietnamica},
  vol.~22, no.~1, 1997.

\bibitem{thi2018}
H.~A. Le~Thi and T.~Pham~Dinh, ``{DC programming and DCA: thirty years of
  developments},'' \emph{Mathematical Programming}, vol. 169, no.~1, 2018.

\bibitem{Auslender1996}
A.~Auslender, ``{Noncoercive Optimization Problems},'' \emph{Mathematics of
  Operations Research}, vol.~21, no.~4, 1996.

\bibitem{banerjee2024}
S.~Banerjee, S.~Peddabomma, R.~Srivastava, and A.~Rajwade, ``A likelihood based
  method for compressive signal recovery under gaussian and saturation noise,''
  \emph{Signal Processing}, vol. 217, 2024.

\bibitem{blake1987}
A.~Blake and A.~Zisserman, \emph{Visual Reconstruction}.\hskip 1em plus 0.5em
  minus 0.4em\relax MIT Press, 1987.

\bibitem{nikolova1998}
M.~Nikolova, ``Estimation of binary images by minimizing convex criteria,'' in
  \emph{IEEE ICIP}, vol.~2, 1998.

\bibitem{nikolova1999}
------, ``{Markovian reconstruction using a GNC approach},'' \emph{IEEE
  Transactions on Image Processing}, vol.~8, no.~9, 1999.

\bibitem{tibshirani1996}
R.~Tibshirani, ``{Regression Shrinkage and Selection via the Lasso},''
  \emph{Journal of the Royal Statistical Society. Series B (Methodological)},
  vol.~58, no.~1, 1996.

\bibitem{feng2020}
Y.~Feng, B.~Ding, H.~Graber, and I.~Selesnick, ``{Transient Artifacts
  Suppression in Time Series via Convex Analysis},'' in \emph{{Signal
  Processing in Medicine and Biology: Emerging Trends in Research and
  Applications}}, I.~Obeid, I.~Selesnick, and J.~Picone, Eds.\hskip 1em plus
  0.5em minus 0.4em\relax Springer, 2020.

\bibitem{heng2025}
Q.~Heng, X.~Liu, and E.~C. Chi, ``{Anderson Accelerated Operator Splitting
  Methods for Convex-nonconvex Regularized Problems},'' \emph{arXiv preprint
  (2502.14269v1)}, 2025.

\bibitem{kitahara2021}
D.~Kitahara, R.~Kato, H.~Kuroda, and A.~Hirabayashi, ``{Multi-Contrast CSMRI
  Using Common Edge Structures with LiGME Model},'' in \emph{EUSIPCO}, 2021.

\bibitem{yata2024}
W.~Yata and I.~Yamada, ``{Imposing Early and Asymptotic Constraints on LiGME
  with Application to Nonconvex Enhancement of Fused Lasso Models},'' in
  \emph{IEEE ICASSP}, 2024, see arXiv:2309.14082v2 for revised version.

\bibitem{kuroda2024}
H.~Kuroda, ``{A Convex-Nonconvex Framework for Enhancing Minimization Induced
  Penalties},'' \emph{arXiv preprint (2407.14819v2)}, 2024.

\bibitem{shoji2025}
S.~Shoji, W.~Yata, K.~Kume, and I.~Yamada, ``{An LiGME Regularizer of
  Designated Isolated Minimizers - An Application to Discrete-Valued Signal
  Estimation},'' \emph{IEICE Transactions on Fundamentals of Electronics,
  Communications and Computer Sciences}, 2025.

\bibitem{groetsch1972}
C.~Groetsch, ``{A note on segmenting Mann iterates},'' \emph{Journal of
  Mathematical Analysis and Applications}, vol.~40, no.~2, 1972.

\bibitem{Goberna2010}
M.~A. Goberna, E.~González, J.~E. Martínez-Legaz, and M.~I. Todorov,
  ``Motzkin decomposition of closed convex sets,'' \emph{Journal of
  Mathematical Analysis and Applications}, vol. 364, no.~1, 2010.

\bibitem{CAaMOTiH}
H.~H. Bauschke and P.~L. Combettes, \emph{Convex Analysis and Monotone Operator
  Theory in Hilbert Spaces}, 2nd~ed.\hskip 1em plus 0.5em minus 0.4em\relax
  Springer, 2017.

\bibitem{rao2007}
K.~R. Rao and P.~Yip, \emph{Discrete Cosine Transform: Algorithms, Advantages,
  Applications}.\hskip 1em plus 0.5em minus 0.4em\relax Academic press, 1990.

\end{thebibliography}

\end{document}